\documentclass[11pt]{article}

        \usepackage{authblk}
        \usepackage[numbers]{natbib}
        \usepackage{keyword}

        \usepackage{preamble}
        \usepackage[
            left=2.8cm,   
            right=2.8cm,  
            top=2.7cm,    
            bottom=2.7cm, 
        ]{geometry}

\begin{document}

            \author[1]{Charly Robinson La Rocca\thanks{Email: charly.robinson.la.rocca@umontreal.ca; Corresponding author.}}
            \author[2]{Jean-François Cordeau}
            \author[1]{Emma Frejinger}

            \affil[1]{Department of Computer Science and Operations Research, Université de Montréal, Canada}

            \affil[2]{Department of Logistics and Operations Management, HEC Montréal, Canada}

\title{Combining supervised learning and local search for the multicommodity capacitated fixed-charge network design problem}
\maketitle

\begin{abstract}
    The multicommodity capacitated fixed-charge network design problem has been extensively studied in the literature due to its wide range of applications. Despite the fact that many sophisticated solution methods exist today, finding high-quality solutions to large-scale instances remains challenging. In this paper, we explore how a data-driven approach can help improve upon the state of the art. By leveraging machine learning models, we attempt to reveal patterns hidden in the data that might be difficult to capture with traditional optimization methods. For scalability, we propose a prediction method where the machine learning model is called at the level of each arc of the graph. We take advantage of off-the-shelf models trained via supervised learning to predict near-optimal solutions. Our experimental results include an algorithm design analysis that compares various integration strategies of predictions within local search algorithms. We benchmark the ML-based approach with respect to the state-of-the-art heuristic for this problem. The findings indicate that our method can outperform the leading heuristic on sets of instances sampled from a uniform distribution. 
\end{abstract}
    

\newpage

\newcommand{\fcn}{MCFND}
\newcommand{\xijk}{x_{ij}^k}
\newcommand{\yij}{y_{ij}}
\newcommand{\ryij}{\tilde{y}_{ij}}
\newcommand{\yb}{\boldsymbol{y}}
\newcommand{\xb}{\boldsymbol{x}}

\newcommand{\ry}{\tilde{\yb}}
\newcommand{\yss}{\yb^{\text{SS}}}
\newcommand{\yls}{\yb^{\text{LS}}}

\newcommand{\yssij}{y_{ij}^{\text{SS}}}
\newcommand{\ybss}{\boldsymbol{y}^{\text{SS}}}

\newcommand{\byij}{\bar{y}_{ij}}
\newcommand{\by}{\bar{\yb}}
\newcommand{\yinc}{\by}

\newcommand{\yp}{\yb'}

\newcommand{\cijk}{c_{ij}^k}
\newcommand{\cij}{c_{ij}}

\newcommand{\fij}{f_{ij}}
\newcommand{\uij}{u_{ij}}
\newcommand{\Xij}{X_{ij}}
\newcommand{\dpij}{\delta^p_{ij}}
\newcommand{\Paths}{P^k}
\newcommand{\path}{p}
\newcommand{\zkp}{z^k_p}

\newcommand{\np}{N_i^{+}}
\newcommand{\nminus}{N_i^{-}}
\newcommand{\Pf}{\mathcal{P}}
\newcommand{\Plin}{\widetilde{\mathcal{P}}}
\newcommand{\pij}{\rho_{ij}}

\newcommand{\binloss}{\mathcal{L}_B}
\newcommand{\ls}{\text{LS}^*}

\newcommand{\wa}{\omega_1}
\newcommand{\wb}{\omega_0}

\newcommand{\mlmodel}{\psi}
\newcommand{\trainedmlmodel}{\hat{\psi}}
\newcommand{\fit}{\text{fit}}
\newcommand{\featurizer}{\phi}
\newcommand{\labels}{\mathcal{Y}}

\newcommand{\lsa}{\mathcal{A}_{\text{LS}}}

\newcommand{\spla}{\mathbb{A}_{\text{SPL}}}
\newcommand{\splx}{\mathcal{X}}

\newcommand{\varname}[1]{\texttt{\detokenize{#1}}}

\newcommand{\refeq}[1]{(\ref{#1})}

\newcommand{\lgbma}{LGBMW1}
\newcommand{\lgbmb}{LGBMW-1}

\section{Introduction}
The multicommodity capacitated fixed-charge network design problem (\fcn) is an influential problem in combinatorial optimization (CO). 
It is an NP-hard problem~\citep{magnanti_network_1984} that has been studied for decades and for which many exact algorithms~\citep{chouman_commodity_2017, crainic_exact_2021} and heuristics~\citep{yaghini_hybrid_2013, yaghini_aco-based_2014,  katayama_mip_2020} have been proposed. The research interest for the \fcn\ is motivated by its general structure which is compatible with applications in logistics, telecommunications, and infrastructure planning~\citep{minoux_networks_1989, quilliot_network_2014}. Recently, extensions to the network design problem have been developed to model congestion~\citep{paraskevopoulos_congested_2016, karimi-mamaghan_hub-and-spoke_2020}, scheduling~\citep{hewitt_new_2023, kidd_relax-and-restrict_2024}, multimodality~\citep{basallo-triana_intermodal_2023, real_multimodal_2021} and uncertainty~\citep{an_two-phase_2016, ghanei_two-stage_2023}. At a high level, the objective of the \fcn\ is to jointly determine the structure of the network and the flows that minimize the sum of both the fixed installation costs and variable routing costs. The fixed costs are related to the activation of arcs. They can, for example, represent the decision to open a facility or a transportation route.  The variable costs model the cost of sending a unit of flow on a given arc of the network. A solution to the \fcn\ must respect the capacity available on each arc and the demand for each commodity. 

Recent advancements in heuristic methods~\citep{gendron_matheuristics_2018, katayama_mip_2020} have enabled the discovery of high-quality solutions to large-scale instances. They often exploit continuous relaxations of \fcn, which are computationally efficient to solve~\citep{crainic_slope_2004, katayama_capacity_2009}. The performance of heuristics is often the result of a combination of many different techniques carefully designed by experts. Typically, improving the performance of a heuristic requires a deep understanding of the problem under study and significant research efforts. This limits the rate of improvement of the state of the art (SOA). Given these circumstances, we explore how a data-driven approach can help improve upon the SOA. Our method feeds a prediction from a machine learning (ML) model to a Local Search (LS) algorithm. The ML model is trained via supervised learning and labels are built using near-optimal solutions. This approach is inspired by the recent successes of ML for CO~\citep{bengio_machine_2021}. 

Our work differentiates itself from the current trend in ML, which is to improve performance via increasingly larger models and datasets.  This trend can be explained by scaling laws which suggest that performance strongly depends on the number of parameters, the dataset size and the amount of computing power~\citep{kaplan_scaling_2020}. Many challenges arise when doing so such as increased cost, carbon emissions impact~\citep{dhar_carbon_2020} and reproducibility concerns~\citep{albertoni_reproducibility_2023}. These issues are well known and many researchers are actively working on sample-efficient methods to obtain accurate predictions with small datasets~\citep{wallingford_neural_2023, robinson_la_rocca_one-shot_2024}. Given that one of our objectives is to remain reproducible, we choose to work with off-the-shelf ML models that are widely accessible. Our contribution lies in the design of a technique for the integration of such models within a heuristic algorithm. In particular, we propose a prediction method where the ML model is called at the level of each arc of the graph. This ensures that the output complexity does not scale exponentially with the size of the instance. This is a key feature of our method that makes it computationally efficient, even for large instances. 


The open research question we attempt to answer is how to best use the ML prediction to help guide a SOA heuristic for the \fcn. We believe that our work is the first to explore this question as most of the literature on this problem focuses on exact methods and human-designed heuristics.  To bridge this gap, this paper makes the following contributions: 
\startblue
\begin{itemize}
    \item A framework for the integration of supervised learning and CO. 
    \item A set of sampling routines that generate informative features for the \fcn. 
    \item A study of the performance of the SOA heuristic with and without predictions of different ML models.
\end{itemize}
\stopblue

The remainder of this article is organized as follows.\startblue Section~\ref{sec:fcn_problem_formulation} presents the problem formulation and related works on the \fcn\ are described in Section~\ref{sec:related_works}.\stopblue We also review some of the literature on supervised learning and the integration of ML in CO. This is followed in Section~\ref{sec:methodology} by our methodology, which includes a general framework for hybrid ML-CO and the sampling routines used to generate informative features. In Section \ref{sec:results}, we quantify our performance with metrics that measure the quality of solutions and the speed at which they are found. Our test set includes the GT dataset~\citep{hewitt_combining_2010} and large instances produced by the Canad generator~\citep{larsen_pseudo-random_2023}. For each instance, we compare the performance of the SOA method~\citep{katayama_mip_2020} augmented by ML predictions against its standalone version. Finally, Section~\ref{sec:conclusion} concludes the paper and discusses future work.

\startblue 
\section{Problem Formulation} \label{sec:fcn_problem_formulation}
\stopblue 

The \fcn\ can be defined on a graph \( G = (N, A) \), where \( N \) is the set 
of nodes and \( A \) is the set of arcs. Each commodity $k \in K$ has a demand $d_k$ that must be routed from its origin $O(k) \in N$ to its destination $D(k) \in N$. The amount of flow on arc $(i,j) \in A$ for commodity $k$ is denoted by $\xijk$ and it incurs a cost per unit of flow $\cijk$. If any amount of flow is routed on arc $(i,j)$, a fixed cost $\fij$ must be paid to use it. This decision is represented by the binary variable~$\yij$.  The objective  is to minimize the sum of both the fixed and variable costs of the arcs used to route all commodities. In terms of constraints, the total flow on arc~$(i,j)$ cannot exceed its capacity~$\uij$. Also, flow conservation must be satisfied at each node $i \in N$. To model this constraint, we define sets of successors $\np$ and predecessors $\nminus$ for each node~$i \in N$. The arc-flow formulation (AF) of $\Pf$ can be stated as the following mixed-integer program (MIP): 
\begin{align}
(\Pf) \quad & \nonumber \\
\min \quad & \sum_{k \in K} \sum_{(i,j) \in A} \cijk \xijk + \sum_{(i,j) \in A} \fij \yij \label{eq:obj} \\
\text{s.t.} \quad & \nonumber \\ 
& \sum_{k \in K} \xijk \leq u_{ij} \yij,  \quad \forall (i,j) \in A, \label{eq:capacity} \\
& \sum_{j \in \np} \xijk - \sum_{j \in \nminus} x_{ji}^k = \begin{cases}
    d^k, & \text{if } i = O(k) \\
    0, & \text{if } i \neq O(k), D(k) \\
    -d^k, & \text{if } i = D(k)
    \end{cases}, \quad i \in N, k \in K, \label{eq:flow_conservation} \\
& \xijk \geq 0, \quad \forall (i,j) \in A, k \in K, \label{eq:nonnegativity} \\
& \yij \in \{0,1\},\quad \forall (i,j) \in A. \label{eq:binary}
\end{align}
The objective function~\refeq{eq:obj} minimizes the sum of flow routing costs and arc opening costs. Constraints~\refeq{eq:capacity} ensure that the total flow on each arc does not exceed its capacity and  bind the decision of opening an arc to the amount of flow routed on it. Constraints~\refeq{eq:flow_conservation} enforce flow conservation at each node. Finally, constraints~\refeq{eq:nonnegativity} and \refeq{eq:binary} define the domains of the decision variables.

\section{Related Works} \label{sec:related_works}

This section describes the body of work related to the \fcn, with an emphasis on heuristics developed to identify near-optimal solutions in a short amount of time.  We then transition to a review of supervised learning, and the gradient boosting models that we use in our experiments. The final part of this section synthesizes recent efforts related to the integration ML and CO.


\subsection{Continuous Relaxations for the \fcn} 
It is common for \fcn\ heuristics to exploit a continuous relaxation of problem $\Pf$~(defined by \refeq{eq:obj} - \refeq{eq:binary}) where the integrality requirements on the binary variables $\yij$ are relaxed. The relaxed problem $\Plin$ is significantly easier than the original problem $\Pf$ and its optimal solution can be used as a tool to explore to solution space efficiently. The main hypothesis is that the solution of $\Plin$ is a reasonable approximation of the optimal solution for $\Pf$. In particular, if the total amount of flow that goes through a given arc in $\Plin$ is zero, the corresponding binary variable is likely to be zero in the optimal solution to the original problem. Inversely, if most of the capacity on the arc is used in $\Plin$, the corresponding binary variable is likely to be one in the optimal solution. This insight is the basis for scaling approaches that iteratively solve the relaxed problem with a modified objective function~\citep{crainic_slope_2004} or set of constraints~\citep{katayama_capacity_2009}. 

 A popular method in this category is Slope Scaling (SS)~\citep{crainic_slope_2004}, which  adjusts the objective function coefficients of each arc to reflect the effective cost of passing flow through it, by considering a combination of both fixed and variable costs. More precisely, the objective function coefficient of an arc for a given commodity becomes the sum of the variable cost and a linearization factor. This factor is the ratio of the fixed cost and the total flow on the arc. The idea is to increase the coefficient of arcs for which the fixed cost is high relative to the flow. One can also set the linearization factor to a large number to force the flow to be zero on certain arcs. Given the optimal solution $\tilde{\xb}$ of $\Plin$, we can construct a feasible solution $\ybss$ for the original problem by rounding up the normalized total flow on each arc. Here, bold formatting is used to denote vectors. The rounding is done according to the following equation:
\begin{align} \label{eq:slope_scaling}
    \yssij = \left\lceil \frac{\sum_{k \in K} \tilde{x}_{ij}^k}{u_{ij}} \right\rceil, \quad \forall (i,j) \in A. 
\end{align}
    
Capacity Scaling (CS)~\citep{katayama_capacity_2009} is a similar approach that iteratively adjusts the capacity of each arc. The method first solves the continuous relaxation of the path-flow formulation (PF) of the \fcn~to obtain a fractional solution $\ry$. In addition to the design variables, the PF formulation incorporates variables $\zkp$ to represent the flows routed through every possible path $p$ for commodity $k$. These variables are defined for every path in the set $\Paths$, which contains all feasible paths for commodity $k$ from the origin $O(k)$ to the destination $D(k)$. The PF can be stated as follows: 
\begin{align}
    \min \quad &  \sum_{(i,j) \in A}  \sum_{k \in K} \cijk \sum_{\path \in \Paths } \dpij \zkp + \sum_{(i,j) \in A} \fij \yij \label{eq:obj:pf} \\
    \text{s.t.} \quad & \nonumber \\ 
    & \sum_{\path \in \Paths} \zkp = d^k,  \quad \forall k \in K, \label{eq:flow:pf} \\
    & \sum_{k \in K} \sum_{\path \in \Paths} \dpij \zkp \leq u_{ij} \yij,  \quad \forall (i,j) \in A, \label{eq:capacity:pf} \\
    & \sum_{\path \in \Paths } \dpij \zkp  \leq  d^k \yij, \quad \forall (i,j) \in A, k \in K, \label{eq:forcing:pf} \\
    & \zkp \geq 0, \quad \forall \path \in \Paths, k \in K, \label{eq:nonnegativity:pf} \\
    & \yij \in \{0,1\},\quad \forall (i,j) \in A. \label{eq:binary:pf}
\end{align}
The binary parameter of $\dpij$ indicates whether or not the arc $(i,j)$ is included in path $\path$. The objective function \refeq{eq:obj:pf} minimizes the total cost similar to the AF formulation. Constraints~\refeq{eq:flow:pf} act as flow conservation constraints for each commodity. Constraints \refeq{eq:capacity:pf} represent the capacity of each arc and \refeq{eq:forcing:pf} provide the forcing constraints. Non-negativity of path-flow variables is enforced by \refeq{eq:nonnegativity:pf}. 

The CS routine iteratively solves the PF formulation~(defined by \refeq{eq:obj:pf}-\refeq{eq:nonnegativity:pf}) and updates the capacity $\uij$ according to the following equation:
\begin{align}
    \uij' = \lambda \uij \ryij + (1-\lambda)  \uij, \quad \forall (i,j) \in A, 
\end{align}
where $\lambda \in [0,1]$ is the so-called capacity smoothing factor that controls the aggressiveness of the scaling. The idea is to increase the capacity of arcs that are used in the solution and decrease the capacity of arcs that are not used. The method is similar to slope scaling in that it aims to force design variables to be zero or one. If $\ryij$ is integral for all arcs, then the solution process stops as the solution is feasible for the original problem. However, in practice, termination criteria are used to stop the process after a certain number of iterations or when the total number of fractional variables is below a given threshold.

The CS procedure in~\cite{katayama_capacity_2009} integrates a column generation technique to improve its computational efficiency. This is required because the total number of paths in the PF is often huge. The method starts with a small set of paths and iteratively adds new paths to the set to improve the solution. The new paths are identified by solving the pricing problem, which considers the reduced cost of each path flow variable. The method is effective because it allows the solution process to focus on a small set of paths that are likely to be used in the optimal solution. At each iteration of the method, the computational effort increases and, therefore, a termination criterion is used to stop the process after a certain number of iterations.

\subsection{Local Search for the \fcn} 
In the previous section, we introduced different scaling approaches that aim to build a feasible solution for $\Pf$  by iteratively solving a relaxed problem. However, they rarely find the optimal solution. In fact, they are most useful as a preprocessing step to warm-start another heuristic. They are designed to explore the solution space around a given solution to find a better one~\citep{shaw_using_1998}. SOA heuristics for the \fcn~often combine scaling and LS to find high-quality solutions. In this section, we describe two examples of LS heuristics, which are used in our experiments. 

The work from~\cite{gendron_matheuristics_2018} presents a matheuristic that combines iterative linear programming and SS to solve the \fcn. At each iteration, the solution $\ry$ of the linear relaxation is used to build index sets $A^0 = \{(i,j)\in A| \ryij = 0\}$ and $A^1 = \{(i,j)\in A| \ryij = 1\}$. In the SS method, arcs in $A^0$ are forced to be zero by fixing the linearization factor to a large number. Inversely, the linearization factor is fixed to 0 for arcs in $A^1$. A feasible solution $\yss$ is then constructed using~\refeq{eq:slope_scaling}. In the next step, the selected variables corresponding to these index sets are fixed at their LP relaxation values in the restricted problem. This problem is then solved using a MIP solver with $\yss$ as a warm-start solution.  A pseudo-cut is then added to the MIP model to remove from the feasible set all solutions that satisfy $\yij = \ryij, (i,j) \in A^0 \cup A^1$: 
\begin{align}
    \sum_{(i,j) \in A^0} \yij + \sum_{(i,j) \in A^1} (1-\yij) \geq 1.
\end{align} 
The method is compared with many other heuristics such as Cycle-based Tabu search~\citep{ghamlouche_cycle-based_2003}, Path Relinking~\citep{ghamlouche_path_2004}, Capacity Scaling~\citep{katayama_capacity_2009} and Simulated Annealing with column generation~\citep{yaghini_hybrid_2013}. The results show that the method is competitive on the Canad instances~\citep{crainic_bundle-based_2001}, which constitute a widely used benchmark for the \fcn. 

More recently, \cite{katayama_mip_2020} proposed to combine CS with a MIP neighbourhood search algorithm, which is the most effective method to produce high-quality solutions for the \fcn~according to our reading of the literature. Hence, our work uses this SOA method as a benchmark to compare the performance of the ML models and, in the remainder of this manuscript, we refer to this method as $\ls$. The heuristic follows a similar strategy to the one described in~\cite{gendron_matheuristics_2018}. However, there are many differences that help improve the performance of the method. First, the CS method is not only used to find a feasible solution, but also to aggressively prune the problem. The set of arcs associated with the design variables that converged to 0 during the CS phase are explicitly deleted from the problem instead of being fixed. Both the continuous and binary variables are removed from the problem, which significantly reduces the memory and time required to solve the problem. Second, $\ls$ adds objective cut-offs to force the MIP solver to find strictly better solutions than the incumbent. Third, the neighbourhood built around the incumbent solution is dynamically updated based on the result of the previous iteration. This neighbourhood is defined by two constraints. The first one removes the incumbent solution $\by$ from the feasible set and forces the solver to find a new solution with at least one additional arc closed: 
\begin{align} \label{eq:kata1}
    \sum_{(i,j) \in A | \byij = 1} \yij \leq \sum_{(i,j) \in A} \byij - 1.
\end{align}
The second constraint defines the size $M$ of the neighbourhood to explore: 
\begin{align} \label{eq:kata2}
    \sum_{(i,j) \in A | \byij = 1} \yij \geq  \sum_{(i,j) \in A} \byij - M.
\end{align}
It states that the number of open arcs in the new solution must be at most $M$ units smaller than the number of open arcs in the incumbent solution. If the solver fails to find a feasible solution within the time limit, the size of the neighbourhood $M$ is decreased by a factor of 2. In the case where the solver proves that no feasible solution exists in the neighbourhood or $M=0$, the algorithm stops. For more details on $\ls$, we refer to the original paper~\citep{katayama_mip_2020}. Our work explores how ML can be used to improve the performance of $\ls$ by predicting near-optimal solutions via supervised learning.

\subsection{Supervised Learning} 
Supervised learning is a branch of ML that aims to learn a mapping between input and output data. In the context of \mkblue{CO}, the input data is a set of features \mkblue{describing} the problem instance and the output data is related to its solution. Often, supervised learning models for binary classification are trained by minimizing the binary cross-entropy loss 
\begin{align} \label{eq:binloss}
    \binloss (y', y) = -\wa y \log(y') - \wb (1-y) \log(1-y'),
\end{align}
where $y'$ is the predicted label and $y $ is the true label. We use the notation $\wa$ and $\wb$ to denote the weights of the positive and negative classes, respectively. \mkblue{These} weights are used to control the contribution of each class to the loss, which is useful when the classes are unbalanced. We cover here the linear binary and gradient boosting models that are used in our experiments.

Linear binary classifiers are widely appreciated for their simplicity, and efficiency in handling linearly separable datasets. They are defined by a set of weights and a bias term, which are multiplied by the input features and added to the result, respectively. The sign of the result determines the predicted label. Once trained with $\binloss$, they produce a hyperplane that separates the input features into two classes~\citep{hastie_linear_2001}. The model's performance is constrained by the fact that it can only represent linear relationships. It limits its ability to capture complex patterns in the data. This can be mitigated by introducing non-linear transformations of the input features, which is the basis for kernel methods~\citep{hastie_kernel_2001}. 

Gradient Boosted Machines (GBM) represent a class ML algorithms that combine multiple learners to form a \mkblue{single} strong predictive model. The essence of GBMs lies in iteratively correcting the residuals of the ensemble, thereby improving the model's accuracy with each step. The foundational work for GBM~\citep{friedman_greedy_2001} introduced the concept of boosting and demonstrated its efficacy in minimizing various loss functions in a stepwise fashion. Empirical results show that boosted tree methods tend to perform exceptionally well across a wide range of tasks when compared to other alternatives such as support vector machines (SVM) or neural networks~\citep{caruana_empirical_2006, borisov_deep_2022_2}. Building upon this foundation, the development of XGBoost~\citep{chen_xgboost_2016} and, more recently, Light GBM (LGBM)~\citep{ke_lightgbm_2017} marked a significant advancement in the field. These implementations of GBMs include various algorithmic enhancements, such as regularized learning objectives, parallelized tree construction, and efficient memory management.

Given the tabular nature of our dataset that relates to the structure of the \fcn, we chose to focus on GBMs as our primary ML model. In particular, we use LGBM, a widely recognized and accessible library, which contributes to the reproducibility of our work.

\subsection{ML for CO}
There is a vast literature on the use of ML techniques to solve CO problems~\citep{dai_learning_2017, bengio_machine_2021, kotary_learning_2021}. Within MIP solvers, ML predictions can support selection strategies for nodes~\citep{khalil_mip-gnn_2022}, cuts~\citep{turner_adaptive_2023} or branches~\citep{khalil_learning_2016, lodi_learning_2017}. In this section, we focus on the use of ML to inform a heuristic that reduces the search space of the problem with the goal of finding high-quality solutions in a short amount of time. A common strategy is fix-and-optimize~\citep{hanafi_mathematical_2017}, where ML is used to identify which subset of variables to fix. \cite{robinson_la_rocca_one-shot_2024} propose an entropy minimization technique to predict, in a one-shot fashion, the optimal solution of generic MIP instances with SOS1 constraints. Given the tripartite graph structure of a MIP instance, Graph Neural Networks (GNN)~\citep{hamilton_graph_2020, liu_machine_2022} are often used to capture patterns in the data and predict solution values for binary variables. The  work of~\cite{ding_accelerating_2020} demonstrates the effectiveness of this strategy on eight  different classes of CO problems. The authors show that a local branching cut~\citep{fischetti_local_2003} of the form
\begin{align}\label{eq:localbranching}
    \sum_{(i,j) \in A} \yij  (1 - \yij') + (1 - \yij) \yij' \leq \beta, 
\end{align} 
which restricts the search space within a distance $\beta$ of the predicted solution $\yp$, can improve the primal gap when compared to SCIP's aggressive heuristic setting. They also compare the precision of their GNN model with XGBoost and show that they have similar performance when both models have access to the LP relaxation solution. This suggests that GBMs can be a reasonable alternative to GNNs even if the data has a graph structure.

Their methodology is applicable to a wide range of problems and, for that reason, their experimental results are compared to a generic MIP solver. For many CO problems, specialized algorithms are available that exploit the structure of the problem to find high-quality solutions. The travelling salesman problem (TSP) is a well-known example of a CO problem that has been studied extensively and for which specialized solvers, \mkblue{such as Concorde~\citep{applegate2006concorde}, exist.}

To show that ML can be competitive with SOA algorithms, \cite{sun_generalization_2021} propose a reduction strategy for TSP instances that uses an SVM to predict which decision variables to remove. They introduce the idea of stochastic sampling of feasible tours to generate statistical features for the SVM. This is effective because one can generate a feasible route of the TSP in polynomial time. We mention this specific work because we use a similar strategy to generate so-called \emph{sampling features}, which we introduce in Section~\ref{sec:methodology}. Using the sampling data, the authors construct statistical measures that aggregate the information for each arc. Their first measure is computed using a weighted sum, where weights are inversely proportional to the rank of the solution.   Mathematically, the rank-based feature $ \phi^r _{ij}$ for arc $(i,j)$ is defined as follows: 
\begin{align}
    \phi^r _{ij} = \sum_{k=1}^{m} \frac{\yij^k}{r_k}, 
\end{align} 
where $r_k$ is the rank of the solution $k$ based on the objective value, and $m$ is the total number of samples. They also developed a measure that uses the Pearson correlation coefficient between each variable and objective values across samples. Their strategy is motivated by the hypothesis that arcs that frequently appear in high-quality solutions are more likely to be part of the optimal solution. Their results show that the SVM model is competitive when compared to other generic problem reduction methods such as Construct, Merge, Solve and Adapt (CSMA)~\citep{blum_construct_2016}. 

In this paper, our main objective is to study how ML predictions can be used to accelerate a SOA heuristic for the \fcn. From our assessment of the literature, this has not been explored before. We believe that our work is relevant because the \fcn~has a wide range of applications, which means that our results can be relevant to many other related problems.  Furthermore, we aim to propose a methodology that is sample efficient and easily reproducible. In the next sections, we describe the methodology followed by the experimental setup used to evaluate its performance.

\section{Methodology} \label{sec:methodology}

We follow a similar approach that of~\cite{sun_generalization_2021}, where an SVM is trained to predict which arcs should be excluded from TSP instances. The core idea is to solve a set of training instances to near-optimality to collect labels for the supervised learning model.  Our methodological contributions come from the adaptations we \mkblue{make} to the original approach to make it suitable for the \fcn, and from the fact that we integrate ML predictions into a SOA heuristic.

An important characteristic of our method is related to the decision to make predictions locally on each individual arc instead of globally on the entire graph~\citep{liu_machine_2022}. Given an instance with $|A|$ arcs, the output space size is thus $2^{|A|}$ because each arc has two possible outcomes. The exponential growth makes global predictions challenging to scale to large instances, while maintaining the ML approach tractable, efficient and accurate. \mkblue{Thanks to local predictions}, we can use a simple binary classifier as our supervised learning model. The main trade-off is that we lose the ability to capture information from the entire graph. This can be mitigated via feature engineering where we extract graph and sampling features to inject global information into the model. 

\startblue 
\textbf{Feature engineering. } Features are always normalized by the corresponding maximum value to ensure that their values are contained in the interval between 0 and 1. We use  aggregation techniques to compress the information in a fixed number of features, which is required by the ML model. We provide a summary of the features used in Table~\ref{tab:fcn:feature_groups} with the number of features in each group and a brief description. In the following paragraphs, we explain in more detail the features derived from the graph and the features generated by the sampling routines.

\begin{table}[h]
\small
\centering
\caption{\mkblue{Summary of feature groups used in ML models}}
\label{tab:fcn:feature_groups}
\begin{tabular}{lccc}
\toprule
{Feature Group} & {Count} & {Description} \\
\midrule
Arc features & 3 & Capacity, variable cost, fixed cost \\
Node features (source/sink) & 10 & Degree, supply, supply sign, max inbound/outbound capacity \\
Fractional solutions & 7 & open/closed frequency, frequency in five bins (histogram) \\
Feasible solutions & 10 & Binary values in top 10 best solutions  \\
\bottomrule
\end{tabular}
\end{table}

\textbf{Graph features. }  Each arc in the \fcn~has three basic characteristics: the capacity, the variable cost and the fixed cost. This first group contains few features and it is not particularly informative about the structure of the instance. For that reason, we also include a second group of features that represent the local topology of the graph. For each arc, we aggregate the information of head and tail nodes. A node can be characterized by its degree and supply, which are related to the  arcs connected to it. We also include the maximum capacity of inbound and outbound arcs.  Note that the exact design choice for graph features is not critical to the performance of the heuristic. As we show in Section~\ref{sec:results}, local graph features are weak predictors of the optimal class. Hence, we direct the feature engineering procedure towards solution sampling techniques.

\textbf{Sampling features. } The idea to sample solutions to generate features is inspired by~\cite{sun_generalization_2021}, where random TSP routes are generated to collect data for the supervised learning model. They construct two measures that aggregate the information using statistical tools.  For our purposes, we include statistical features in a similar fashion. During the solution sampling process, we collect both fractional and feasible solutions. Fractional solutions are easy to generate and we find a large number of them. We compress that information with a histogram method that computes the frequency of fractional values in each of five bins between 0 and 1. The idea is to capture the distribution of fractional values in the solution space. We also include how frequently each arc is fully open or closed in the fractional solutions. Our preliminary experiments used standard aggregations (average, min, max) of the fractional values, but we found that the histogram technique was marginally better. Simple features were surprisingly effective in our preliminary experiments. We find significantly fewer feasible solutions than fractional ones during the sampling process. Thus, it is not necessary to aggregate them. We include as-is the 10 best feasible solutions as binary features.  We use padding with zeros if the number of feasible solutions is smaller than 10.
\stopblue  



\textbf{Sampling routines. } 
Given that, for the \fcn, there is not a well established way to build feasible solutions in polynomial time akin to the TSP, we propose sampling routines to generate solutions. They take advantage of a continuous relaxation of the problem, which is computationally efficient to solve.  They are designed to generate a large and diverse set of solutions in a short amount of time. In particular, we propose Randomized Slope Scaling (RSS) and Local Search Fast Sampling (LSFS), which use SS and CS, respectively. We utilize two different sampling routines to diversify and reduce the bias in the set of samples.

The RSS routine is based on~\cite{gendron_matheuristics_2018},  where a MIP solver and SS are called sequentially. We found that this heuristic is not competitive with $\ls$ in terms of solution quality. However, it is effective as a sampling routine because it can produce feasible solutions relatively quickly. We modified the original implementation in \mkblue{two} different ways to improve its ability to explore the solution space. First, we pass a random value for both the \varname{Seed} and \varname{Heuristics} parameters of Gurobi at each iteration. Second, we set the objective coefficient to a large number for 1\% of the variables chosen at random. 

\startblue 
The LSFS routine uses the same implementation as $\ls$ with a few parameter changes. The adjustments aim to reduce the computational time associated with CS, which scales with the number of columns generated. The specific choice of parameters is provided in Section~\ref{sec:results}.  \stopblue 

\textbf{Training. } The training process is straightforward because we use off-the-shelf ML libraries that have high-level interfaces available. In a supervised setting, we provide the ML model $\mlmodel$ with a set of features and their corresponding labels. The featurizer $\featurizer$ is a function that takes the set of solutions $\splx$ and returns a set of features according to the feature engineering process described above. The set of binary labels $\labels$ are built by merging the top three best known solutions instead of just the best one, i.e., all arcs that are open in at least one of the solutions are labelled as open. Both the features and the labels are passed to the library-specific $\fit$ routine to generate a trained ML model $\trainedmlmodel$ as follows:
\begin{align} \label{eq:fit}
    \trainedmlmodel \gets \fit\left( \mlmodel, \featurizer(\splx), \labels \right). 
\end{align}

\textbf{Inference. } Algorithm~\ref{alg:sls} describes the Supervised Local Search (SLS) algorithm, which is a high-level description of our inference pipeline. It is designed to be flexible, so it can be adapted to different CO problems. It takes as input the problem $\Pf$, the trained ML model $\trainedmlmodel$, the local search algorithm $\lsa$, and the set of sampling routines $\spla$. In the first phase, the algorithm parallelizes the sampling process using the routines in $\spla$. The \mkblue{collected solutions} are accumulated in a set $\splx$, which is then used to get an incumbent solution $\yinc$ that can warm-start the solution process. We get an estimate of near-optimal solutions $\yp$ as the output of the ML model with $\featurizer(\splx)$ as its input. The local search algorithm $\lsa$ is then called to attempt to find a high-quality solution $\yls$ of the problem $\Pf$ . Note that Algorithm~\ref{alg:sls} does not explicitly describe  how the ML predictions are integrated in the solution process. This is an implementation detail of the algorithm $\lsa$, and Section~\ref{sec:algo_design} provides experimental results for \mkblue{specific} strategies. 



\begin{algorithm}
    \caption{Supervised Local Search (SLS) algorithm}\label{alg:sls}
    \begin{algorithmic}
    \Procedure{SLS}{$\Pf, \trainedmlmodel, \lsa, \spla$}
    \State $\splx = \{\}$ \Comment{Initialize the set of solutions} 
    \For{(In parallel)~$\alpha \in \spla$}
        \State $\splx \gets \text{sample}(\Pf, \alpha)$ \Comment{Generate features using the sampling routine}
    \EndFor

    \State $\by \gets \text{get\_best}(\splx)$ \Comment{Get best solution from the set}
    \State $\yp \gets \trainedmlmodel(\featurizer(\splx))$ \Comment{Forward pass of the ML model}

    \State $\yls \gets \text{solve}(\Pf, \lsa, \by, \yp)$ \Comment{Solve using the local search algorithm}
    \State \textbf{return} $\yls$
    \EndProcedure
    \end{algorithmic}
\end{algorithm}

\textbf{Problem reduction. }  The problem reduction step is a key part of the $\ls$ algorithm and one of the main reasons \mkblue{for its effectiveness}. At every CS iteration, each arc with a fractional value below a certain threshold is removed from the problem. This significantly reduces the memory and time required to solve the problem. However, it introduces a challenge related to the integration of \mkblue{predicted} solutions. In particular, we found that both the incumbent $\yinc$ and the ML prediction $\yp$ are rarely feasible in the reduced problem. Furthermore, $\ls$ already incorporates local cuts~\refeq{eq:kata1} and~\refeq{eq:kata2}, which are more restrictive than \refeq{eq:localbranching}. For these reasons, the cut~\refeq{eq:localbranching} does not improve the performance of $\ls$. To circumvent the misalignment between the reduced problem and the \mkblue{predicted} solutions, we decided to replace the original CS-based reduction by a routine that uses ML to predict which arcs to remove, akin to~\cite{sun_generalization_2021}. This is motivated by the fact that the CS routine can take up to 1 hour to execute. Also, we can make the required adjustment to make sure the incumbent and the ML prediction is feasible in the reduced problem. 

\section{Experimental Results} \label{sec:results}

This section presents the results of our experiments. We first describe the experimental setup and the different datasets used to test our method. We then present some intermediate results on the extended runs and the sampling routines, which give us the necessary information to train the ML model. We include results on the tuning procedures for $\ls$ and our ML models. This is followed by an algorithm design analysis where we test different ways to integrate the ML prediction in an LS algorithm. Finally, we explain the main results of our work, which is the comparison of the performance of $\ls$ with and without ML predictions.

\textbf{Metrics. } The metrics we use to evaluate the performance of the heuristics are the primal gap (PG) and the primal integral (PI). The primal gap is the relative difference between the objective function values of the best solution found by the heuristic and of the best known solution. The primal integral is calculated using the area under the curve of the primal gap over time. The primal gap is a measure of the quality of the solutions found, and the primal integral is a measure of the speed at which the heuristic finds high-quality solutions. Each table contains the mean, the shifted geometric mean and quantiles of the corresponding metric.\startblue The geometric mean appears with a shift of one in all tables and we refer to it as geomean in the remaining of this paper.   We characterize the performance of the ML models with the balanced accuracy. It is the average of recall obtained on each class, which helps address the unbalanced class distribution. \stopblue

\textbf{Experimental settings. } Our software stack is written in Julia and we use the JuMP package~\citep{lubin_jump_2023} to interface with the Gurobi Optimizer version 10.1~\citep{gurobi}.  We run each \mkblue{experiment for two hours} on a single thread of an Intel Gold 6148 Skylake (2.4 GHz) processor and 35 GiB of memory. We use the original implementation of $\ls$, which was written in C++ by~\cite{katayama_mip_2020}. \startblue We use the same parameters as the original codebase, except when we mention otherwise. Both sampling routines (RSS, LSFS) have a total time budget of 5 minutes, and a time limit of 20 seconds per iteration. We adjust the number of columns generated in LSFS and we study the impact of that parameter in Section~\ref{sec:extended_runs}. Another important parameter that we consider is the capacity smoothing factor $\lambda$, and we describe the related tuning in Section~\ref{sec:tuning}. This section also includes a discussion about the choice of the ML model and its parameters. 
\stopblue

\textbf{Datasets. } We experiment with sets of instances with different characteristics. The first set contains the 24 GT instances~\citep{hewitt_combining_2010}. They are respectively subdivided into GTFL and GTFT, which stand for instances with loose and tight capacities. We chose to work with these instances because they are challenging and they are the ones used in the original work~\citep{katayama_mip_2020}. The second set contains instances built by a modernized version of the generator used to produce the Canad dataset~\citep{larsen_pseudo-random_2023}. Among others, the generator has parameters to control the size of the instance (number of arcs and commodities), and its topology. We set the instance size to be similar to GT (between 2,000 and 3,000 arcs, and around 200 commodities) to minimize the potential performance impact on $\ls$ due to calibration issues. The sets A and B use the grid topology, whereas C and D use the circular one. The sets A and C are exclusively used for training purposes and the sets B and D are used to test the method. Instances in test sets are about 30\% larger than those in the training sets, to ensure that the method can generalize. Each generated set contains 30 instances, for a total of 144 instances when we include the GT instances.

\subsection{Extended Runs and Sampling Routines} \label{sec:extended_runs}

This section covers how to collect the required data to train the ML model. Features and labels are generated using the sampling routines and extended runs, respectively. Table~\ref{tab:opt_gap_12h} contains the optimality gaps given by the Gurobi default configuration (grb) after 12-hour runs. It gives us an estimate of the relative difficulty of each dataset. We notice how GT instances are meaningfully more challenging than the generated ones. Dataset A instances are the easiest with an average optimality gap of 1.3\%, and dataset GTFL instances are the hardest with an average gap of 27.37\%. The circular typology is marginally more difficult than the grid one. Test instances are slightly more challenging because they are larger.  

Table~\ref{tab:primal_gap_12h} includes primal gaps for Gurobi and $\ls$ on each dataset. The main takeaway is the relative effectiveness of $\ls$ when specifically applied to GT instances. Gurobi has an average gap of 24.84\% on GTFL, and $\ls$ has an average of 0\%. This means the SOA heuristic systematically finds the best solution. The story is different for easier instances where Gurobi tends to outperform $\ls$.  

To generate the features, we run each sampling routine for 5 minutes and collect all the solutions found during the solution process. The results related to sampling routines on GT instances are summarized in Tables~\ref{tab:primal_gap_samplingGT} and \ref{tab:cb_sol_count_samplingGT}. The first table presents the quality of the best solution found by each routine, and the second contains the solution count. Gurobi (GRBP-5) is included as a reference. The LSFS routine is executed with three different column generation values: 3, 10 and 50. We note that LSFS-3 is the most effective routine in terms of solution quality with an average gap of 13.45\%. An increase in the number of columns generated does not help because of the limited time budget. The RSS routine's main advantage is its ability to generate a large number of solutions. It can generate more than 10 times more solutions than LSFS in the same period of time. We found similar results on the generated instances.

\begin{table}[h]
\centering
\caption{Optimality gap (\%) using Gurobi for 12-hour runs}
\label{tab:opt_gap_12h}
\begin{tabular}{lcccccc}
\toprule
{Dataset} & \multicolumn{3}{c}{Quantiles} & {Mean} & {Geomean} \\
{} & {0.1} & {0.5} & {0.9} & {} & {} \\
\midrule
GTFT & 16.58 & 26.36 & 39.92 & 27.37 & 26.88 \\
GTFL & 12.04 & 30.31 & 44.69 & 29.10 & 27.15 \\
A & 0.88 & 1.35 & 1.62 & 1.30 & 2.28 \\
B & 1.85 & 2.46 & 2.90 & 2.40 & 3.37 \\
C & 1.01 & 1.57 & 2.20 & 1.64 & 2.60 \\
D & 2.61 & 3.85 & 4.84 & 3.82 & 4.74 \\
\bottomrule
\end{tabular}
\end{table}

\begin{table}[h]
\centering
\caption{Primal gap (\%) for 12-hour runs with baseline algorithms}
\label{tab:primal_gap_12h}
\begin{tabular}{lcccccc}
\toprule
{Scenario} & \multicolumn{3}{c}{Quantiles} & {Mean} & {Geomean} \\
{} & {0.1} & {0.5} & {0.9} & {} & {} \\
\midrule
grb-GTFT & 0.18 & 13.34 & 31.05 & 14.70 & 9.09 \\
grb-GTFL & 2.51 & 26.33 & 44.84 & 24.84 & 18.45 \\
grb-A & 0.00 & 0.00 & 0.11 & 0.04 & 1.04 \\
grb-B & 0.00 & 0.35 & 0.62 & 0.35 & 1.32 \\
grb-C & 0.00 & 0.00 & 0.13 & 0.02 & 1.02 \\
grb-D & 0.16 & 0.64 & 1.08 & 0.67 & 1.63 \\
LS*-GTFT & 0.00 & 0.00 & 0.91 & 0.23 & 1.15 \\
LS*-GTFL & 0.00 & 0.00 & 0.00 & 0.00 & 1.00 \\
LS*-A & 0.00 & 0.05 & 0.61 & 0.19 & 1.16 \\
LS*-B & 0.00 & 0.00 & 0.51 & 0.12 & 1.10 \\
LS*-C & 0.00 & 0.09 & 0.53 & 0.21 & 1.18 \\
LS*-D & 0.00 & 0.00 & 0.00 & 0.01 & 1.01 \\
\bottomrule
\end{tabular}
\end{table}

\begin{table}[h]
\centering
\caption{Primal gap (\%) for sampling algorithms on dataset GT}
\label{tab:primal_gap_samplingGT}
\begin{tabular}{lcccccc}
\toprule
{Scenario} & \multicolumn{3}{c}{Quantiles} & {Mean} & {Geomean} \\
{} & {0.1} & {0.5} & {0.9} & {} & {} \\
\midrule
LSFS-3 & 0.00 & 1.90 & 35.95 & 13.45 & 4.85 \\
LSFS-10 & 0.00 & 15.17 & 93.13 & 41.60 & 17.84 \\
LSFS-50 & 0.00 & 92.10 & 120.22 & 83.89 & 31.66 \\
RSS & 2.09 & 21.83 & 35.65 & 21.80 & 17.36 \\
GRBP-5 & 12.12 & 26.96 & 40.91 & 25.96 & 24.40 \\
\bottomrule
\end{tabular}
\end{table}

\begin{table}[h]
\centering
\caption{Solution count for sampling algorithms on dataset GT}
\label{tab:cb_sol_count_samplingGT}
\begin{tabular}{lcccccc}
\toprule
{Scenario} & \multicolumn{3}{c}{Quantiles} & {Mean} & {Geomean} \\
{} & {0.1} & {0.5} & {0.9} & {} & {} \\
\midrule
LSFS-3 & 4.30 & 7.50 & 13.70 & 8.29 & 8.69 \\
LSFS-10 & 3.00 & 4.50 & 8.70 & 5.54 & 6.13 \\
LSFS-50 & 2.00 & 2.00 & 6.70 & 3.29 & 3.83 \\
RSS & 19.70 & 53.50 & 206.00 & 104.60 & 66.44 \\
GRBP-5 & 2.00 & 4.00 & 5.00 & 3.71 & 4.53 \\
\bottomrule
\end{tabular}
\end{table}

\subsection{Tuning $\ls$ and the ML Model} \label{sec:tuning}
The original implementation of $\ls$ was tuned on GT instances, so we do not change its default parameters for that dataset. We do, however, tune the capacity smoothing factor $\lambda$ on the newly generated instances. To do so, we solve each instance in the training sets for 30 minutes with different values of $\lambda$. The results for dataset A are presented in Table~\ref{tab:primal_gap_tune_A}. According to the average primal gap, the best value for $\lambda$ is 0.05. It is also the best setting for dataset C. We chose to tune $\lambda$ because it is the parameter that is explicitly calibrated in the original work. Since generated instances have a similar size than GT instances, we expect that the remaining parameters should be adequate.

The ML models are also tuned on the training sets A and C, and we study their performance  with the balanced accuracy metric. Table~\ref{tab:ml_tune_two_cols} reveals the effect of features sets on a selection of ML models. We test different combinations of basic (B), node graph (N) and sampling (S) features.  The basic features are specific to the given arc, whereas node features contain aggregated information about the graph neighbourhood. The sampling features come from the solutions generated by the sampling routines. The linear binary classifier (Lin) is compared to four LGBM models with different values for \varname{max_depth}, which controls its capacity. The main takeaway from this table is the importance of the sampling features. The average validation accuracy increases from 67\% to 92\% when sampling features are included. Note that BNS and BS have a similar validation accuracy, but BNS has a higher training accuracy. This effect is more noticeable as the model capacity increases. This is a sign of overfitting. It suggests that the graph features do not help the model to generalize. 

Note that we considered other parameters of LGBM such as  \varname{min_data_in_leaf}, \varname{learning_rate} and \varname{num_iterations}, but they did not improve the performance. We noticed that LGBM is resilient to overfitting when \varname{min_data_in_leaf} is set to its default value of 20. When we decrease it to 10 or 3, the training accuracy increases, but the validation accuracy decreases. Furthermore, we should mention that LGBM will print a warning ``No further splits with positive gain'', when its depth is higher than or equal to 7. This reinforces the idea that many features are not useful or they are redundant.  LGBM is able to find high-quality splits in the feature space early in the training process, but it struggles to find meaningful splits as the depth increases. Since both the linear model and LGBM have similar performance, the choice of the model is not critical. We choose to use LGBM because it is more flexible and offers more options to manage the class unbalance.

A deeper look into the raw data shows that misclassifications seem to be caused by a weak point in the sampling solutions. In particular, around 10\% of open arcs in the best solution are always closed in sampling solutions. Sampling routines explore the solution space around the optimal solution of continuous relaxations, which is not always near the optimal solution of the original problem. When this is the case, the ML model has to rely on graph features to predict that an arc should be open. It must not be overlooked that the evidence in Table~\ref{tab:ml_tune_two_cols} points to the limited predictive power of graph features. This highlights a limitation of our ML-based approach. The model cannot accurately predict that an arc is open unless there is some evidence of that in the sampling solutions. This is a difficult problem to solve because sampling routines cannot reveal the optimal solution with certainty. 

\textbf{Bias adjustments. } 
An important challenge that emerges when using solutions as labels is the unbalanced nature of the dataset. In a solution of any given \fcn\ instance, it is rarely the case that the number of open arcs is equal to the number of closed arcs. This is a problem because the ML model can learn to predict the majority class and still achieve a high accuracy. To mitigate this issue, one might consider using a resampling technique.  We experimented with a random undersampling technique, but it did not improve performance when the flag \varname{is_unbalance} is set to true in LGBM. The library can handle unbalanced datasets by itself to some extent. The default binary loss, however, is designed to minimize misclassifications with an equal weight for both classes. We suspect that this is inadequate because false negatives are not necessarily as costly as false positives. To address this, we test different weights in the loss~\refeq{eq:binloss} to skew the model in favor of one class over the other.  This incentivizes the model to classify ambiguous cases as positive or negative depending on the bias. We study its effect with two models: \lgbma\ is biased towards the positive class ($\wa = 0.75$, $\wb=0.25$), whereas \lgbmb\ is biased towards the negative class ($\wb = 0.75$, $\wa=0.25$). We quantify the bias via false positive and false negative rates in Tables~\ref{tab:false_negative_rate} and \ref{tab:false_positive_rate}, respectively. Note that models are trained to predict the optimal reduction (see Section~\ref{sec:algo_design} for the justification). Hence, the positive label corresponds to removing an arc from the graph. As expected, \lgbma\ has the lowest false negative rate at 3\% and  \lgbmb\ has the lowest false positive rate at 5\%. We demonstrate in Section~\ref{sec:main_results} that including the bias adjustments can improve performance.




\begin{table}[h]
\centering
\caption{Primal gap (\%) for tuning $\lambda$ on dataset A}
\label{tab:primal_gap_tune_A}
\begin{tabular}{lcccccc}
\toprule
{Scenario} & \multicolumn{3}{c}{Quantiles} & {Mean} & {Geomean} \\
{} & {0.1} & {0.5} & {0.9} & {} & {} \\
\midrule
LS*-0.05 & 0.00 & 0.00 & 0.00 & 0.00 & 1.00 \\
LS*-0.15 & 0.02 & 0.05 & 0.13 & 0.07 & 1.07 \\
LS*-0.25 & 0.08 & 0.14 & 0.23 & 0.15 & 1.15 \\
LS*-0.35 & 0.07 & 0.19 & 0.27 & 0.18 & 1.17 \\
LS*-0.45 & 0.14 & 0.24 & 0.35 & 0.26 & 1.25 \\
\bottomrule
\end{tabular}
\end{table}

\begin{table}[h]
\centering
\caption{Train and validation balanced accuracy for different ML models and sets of features (bold underscore for best values)}
\label{tab:ml_tune_two_cols}
\begin{tabular}{lccccccc}
\toprule
{Model} & \multicolumn{2}{c}{BN} & \multicolumn{2}{c}{BNS} & \multicolumn{2}{c}{BS} \\
{} & {Train} & {Val} & {Train} & {Val} & {Train} & {Val} \\
\midrule
Lin & 0.65 & 0.65 & \underline{\textbf{0.92}} & \underline{\textbf{0.92}} & \underline{\textbf{0.92}} & \underline{\textbf{0.92}} \\
LGBM-D3 & 0.66 & 0.66 & \underline{\textbf{0.92}} & 0.91 & \underline{\textbf{0.92}} & \underline{\textbf{0.92}} \\
LGBM-D7 & 0.70 & 0.68 & \underline{\textbf{0.93}} & \underline{\textbf{0.92}} & \underline{\textbf{0.93}} & \underline{\textbf{0.92}} \\
LGBM-D11 & 0.76 & 0.68 & \underline{\textbf{0.96}} & 0.91 & 0.94 & \underline{\textbf{0.92}} \\
LGBM-D15 & 0.81 & 0.66 & \underline{\textbf{0.96}} & \underline{\textbf{0.92}} & 0.95 & \underline{\textbf{0.92}} \\
\textbf{num best out of 5} & 0 & 0 & 5 & 3 & 3 & 5 \\
\textbf{mean} & 0.72 & 0.67 & 0.94 & 0.92 & 0.93 & 0.92 \\
\bottomrule
\end{tabular}
\end{table}

\begin{table}[h]
\centering
\caption{False negative rate for different ML models}
\label{tab:false_negative_rate}
\begin{tabular}{lcccccc}
\toprule
{Scenario} & \multicolumn{3}{c}{Quantiles} & {Mean} & {Geomean} \\
{} & {0.1} & {0.5} & {0.9} & {} & {} \\
\midrule
Lin & 0.02 & 0.07 & 0.42 & 0.12 & 1.12 \\
LGBM & 0.04 & 0.10 & 0.39 & 0.16 & 1.15 \\
LGBMW1 & 0.01 & 0.02 & 0.05 & 0.03 & 1.03 \\
LGBMW-1 & 0.11 & 0.21 & 0.35 & 0.23 & 1.22 \\
\bottomrule
\end{tabular}
\end{table}

\begin{table}[h]
\centering
\caption{False positive rate for different ML models}
\label{tab:false_positive_rate}
\begin{tabular}{lcccccc}
\toprule
{Scenario} & \multicolumn{3}{c}{Quantiles} & {Mean} & {Geomean} \\
{} & {0.1} & {0.5} & {0.9} & {} & {} \\
\midrule
Lin & 0.06 & 0.11 & 0.48 & 0.19 & 1.18 \\
LGBM & 0.06 & 0.11 & 0.30 & 0.15 & 1.14 \\
LGBMW1 & 0.17 & 0.27 & 0.45 & 0.29 & 1.29 \\
LGBMW-1 & 0.00 & 0.03 & 0.13 & 0.05 & 1.05 \\
\bottomrule
\end{tabular}
\end{table}

\subsection{Algorithm Design Analysis} \label{sec:algo_design}
ML-OR hybrid approaches \mkblue{often} come with a unique challenge, which is related to the integration of a noisy prediction in a \mkblue{search} algorithm. To obtain good performance, the algorithm must be designed in a way that is resilient to a relatively high level of noise.  Related works~\citep{ding_accelerating_2020} use the local branching cut~\refeq{eq:localbranching} to guide the search. We perceive it as a good choice in this context because it is a form of soft variable fixing. For that reason, our first experimentations revolve around the integration of that cut into the original implementation of $\ls$ to guide the search around the ML prediction. 

In this section, we introduce an algorithm design analysis where we test different ways to integrate the ML predictions in an LS algorithm. To simulate the effect of the ML model, we use our database of solutions (DB) to get the correct label and then we add noise to it. The noise is a parameter of the simulation, and it is the probability that the label is flipped. We test four noise values between 0 and 0.15, which are close to the misclassification rate on the validation set. The three algorithms we test are: 
\begin{itemize}
    \item \textbf{LBH}: A local branching cut~\refeq{eq:localbranching} is added relative to the ML prediction. The distance parameter is set such that 80\% of arcs must respect the ML solution. We also set the variable parameters \varname{VarHintVal} in Gurobi to the ML prediction.
    \item \textbf{LSWSH}: Similar to LBH, but the cut~\refeq{eq:localbranching} and the hints are added to $\ls$ instead of Gurobi alone. We use the same distance parameter as LBH. 
    \item \textbf{LSR}: This algorithm reduces the problem instead of guiding the search. Arcs are removed from the graph according to the ML prediction. However, if the reduction is applied as-is, the reduced problem is unlikely to be feasible. To guarantee feasibility, we add a repair routine that keeps arcs that are part of the best sampling solution. LSR uses the $\ls$ implementation without the CS routine. 
\end{itemize}
We compare each algorithm on the GT dataset in Table~\ref{tab:primal_gap_db_GT}. $\ls$-3 is a version of $\ls$ with only 3 columns generated, which significantly reduces the runtime of the CS routine. As expected, the best performing algorithms use DB with no noise. LSWSH, which adds a cut in $\ls$, marginally benefits from the provided solution. This can be explained by an alignment issue. The best known solution is not necessarily feasible in the reduced problem produced by CS. LBH has most of the wins when the noise is null, but its performance decreases rapidly as the noise increases. The reason LBH performs so well with no noise is because \varname{VarHintVal} is set to a high-quality solution that is feasible. LSR is the most resilient algorithm to noise. It is the only algorithm that is able to outperform $\ls$ with a noise of 5 or 10\%. This is a promising result because it suggests that the ML model can be used to reduce the problem size without significantly affecting the quality of the solutions as long as the misclassification rate is not above 10\%. \startblue Given that LSR is the most promising method, we select it for further experimentations and highlight it in bold in the remaining tables of this paper. \stopblue

\begin{table}[h]
\centering
\caption{\mkblue{Primal gap (\%) for algorithms using DB on dataset GT for 2-hour runs}}
\label{tab:primal_gap_db_GT}
\begin{tabular}{lccccccc}
\toprule
{Scenario} & \multicolumn{3}{c}{Quantiles} & {Mean} & {Geomean} & {Wins} \\
{} & {0.1} & {0.5} & {0.9} & {} & {} & {} \\
\midrule
LS* & 0.39 & 1.36 & 16.59 & 4.49 & 3.25 & 0 \\
LS*-3 & 5.65 & 11.83 & 30.01 & 14.46 & 12.95 & 0 \\
LSWSH-DB-0.0 & 0.19 & 0.70 & 42.94 & 17.32 & 3.21 & 0 \\
\textbf{LSR-DB-0.0} & 0.00 & 0.00 & 0.17 & 0.08 & 1.07 & 7 \\
LBH-0.8-DB-0.0 & 0.00 & 0.00 & 0.00 & 0.02 & 1.02 & 10 \\
LSWSH-DB-0.05 & 0.19 & 2.64 & 129.77 & 28.07 & 6.34 & 1 \\
\textbf{LSR-DB-0.05} & 0.45 & 2.48 & 6.03 & 3.15 & 3.53 & 1 \\
LBH-0.8-DB-0.05 & 1.61 & 8.24 & 14.42 & 8.33 & 7.45 & 0 \\
LSWSH-DB-0.1 & 0.38 & 1.48 & 132.03 & 32.58 & 6.00 & 0 \\
\textbf{LSR-DB-0.1} & 0.68 & 2.13 & 5.54 & 2.51 & 3.11 & 1 \\
LBH-0.8-DB-0.1 & 1.26 & 11.51 & 32.92 & 15.77 & 10.96 & 0 \\
LSWSH-DB-0.15 & 0.38 & 6.13 & 134.30 & 30.49 & 8.38 & 0 \\
\textbf{LSR-DB-0.15} & 1.44 & 6.27 & 10.31 & 6.20 & 6.22 & 0 \\
LBH-0.8-DB-0.15 & 1.63 & 11.51 & 32.92 & 16.46 & 11.71 & 0 \\
\bottomrule
\end{tabular}
\end{table}

\subsection{Main Results} \label{sec:main_results}
In this section, we assemble the pieces and compare the performance of $\ls$ with LSR where the ML prediction is used to reduce the problem size. We begin our analysis with the results on GT instances. Tables~\ref{tab:primal_gap_ml_GT_2-hour} and \ref{tab:primal_integral_ml_GT_2-hour} contain the primal gaps and primal integrals, respectively. We select a 2-hour time limit because at least one hour is required to run CS on some GT instances. For reference, we also include Gurobi with \varname{MIPFocus} set to 1 (grb-f1), and LB with a cut around the incumbent. We experiment with three ML models for LSR: LGBM, \lgbma\ and \lgbmb. In terms of primal gap, $\ls$ dominates with 16 wins out of 24 and 1.94 geomean. It is followed by LSR-\lgbmb\ with 3 wins and 4.94 geomean. The primal integral, however, reveals a different story. $\ls$-3 finds high-quality solutions early with 12 wins and a geomean of 5.5. LSR-\lgbmb\ is the second best with 7 wins and a geomean of 7.13. The best primal integral is achieved by $\ls$-3, but at the cost of a higher primal gap. The results reveal that LSR-\lgbmb\ has a good balance between quality and speed. The outperformance of $\ls$ can be justified by two factors. First, $\ls$ was extensively tuned on GT instances.  Our experiments suggest that the CS routine is not as effective on the other instances that we generated. Second, each instance in GT has a unique structure, where the number of arcs and commodities varies significantly from one instance to the other. This is particularly challenging for ML methods, which are known to struggle when the underlying distribution of the data is not consistent.


We now turn our attention to the generated instances, which are more representative of the target application for this method. ML-CO hybrid approaches seek to perform well on instances sampled from a relatively narrow distribution, which is the case for the generated instances.  Tables~\ref{tab:primal_gap_ml_B_2-hour}-\ref{tab:primal_integral_ml_D_2-hour} contain the metrics of interest on test sets B and D. In both cases, LSR-\lgbmb\ outperforms $\ls$ in terms of quality and speed. Its average primal gap is 0.06\% and 0.11\%, which correspond to 14 and 16 wins, for datasets B and D, respectively. $\ls$ is the second-best algorithm with respect to this metric. In terms of primal integral,  LSR-\lgbmb\ is also the best-performing algorithm with 32 wins (17+15) in total for both datasets. \startblue  In Tables~\ref{tab:fcn:gap_over_time_GT}, \ref{tab:fcn:gap_over_time_B} and \ref{tab:fcn:gap_over_time_D}, we present the shifted geomean of the primal gap over time for datasets GT, B and D, respectively. For challenging instances, $\ls$ does not find feasible solutions early because the CS routine is time-consuming. On average, for dataset D, LSR-\lgbmb\ can find solutions of similar quality to $\ls$ roughly one hour earlier.  \stopblue Furthermore, we remark that LB systematically outperforms grb-f1. Hence, it is a baseline that is relatively easy to implement and worth considering. Finally, we ought to mention that the performance of LSR-LGBM is competitive but it is consistently outperformed by LSR-\lgbmb. This is a sign that the bias adjustments are beneficial.

\startblue 
We draw attention to the trade-off that the ML component brings compared to classical heuristic methods. Learning an instance reduction strategy requires high-quality solutions and a training procedure. In our case, the computational cost of the latter is negligible compared to the former. Our ML model takes less than a minute to train on a standard CPU, which is small compared to the time required to solve any given instance. However, the cost associated with the labelling can be a source of concern. For example, in our experimental setup, we spend 12 hours per instance, which amounts to 1,728 hours of computation for the 144 instances. If the current performance obtained using classical heuristics is satisfactory, then there is no need to use ML. However, if the performance is inadequate, then ML can be used to improve it. The cost versus benefit analysis will shift in different contexts where instances are repetitively solved, in time-critical applications, or when computational resources are highly available. We believe that for many real-world scenarios the one-time investment is justified by the performance improvement. \stopblue 

\textbf{Key insights.} $\ls$ dominates on GT instances, as it excels in terms of primal gap but its slow CS routine negatively impacts its primal integral. LSR replaces CS with ML predictions which allows it to find high-quality solutions faster. Its primal integral is an order of magnitude better than $\ls$ on GT. On generated instances, LSR outperforms for both metrics, particularly when using the ML model with bias adjustments. 

\begin{table}[h]
\centering
\caption{\mkblue{Primal gap (\%) with ML models for 2-hour runs on dataset GT}}
\label{tab:primal_gap_ml_GT_2-hour}
\begin{tabular}{lccccccc}
\toprule
{Scenario} & \multicolumn{3}{c}{Quantiles} & {Mean} & {Geomean} & {Wins} \\
{} & {0.1} & {0.5} & {0.9} & {} & {} & {} \\
\midrule
LB-0.8 & 3.34 & 24.87 & 41.15 & 23.83 & 18.85 & 0 \\
LS* & 0.00 & 0.00 & 10.16 & 8.24 & 1.94 & 16 \\
LS*-3 & 0.00 & 9.71 & 33.19 & 13.57 & 8.30 & 5 \\
\textbf{LSR-LGBM} & 2.78 & 8.19 & 21.95 & 10.76 & 9.34 & 0 \\
\textbf{LSR-LGBMW-1} & 0.06 & 5.91 & 9.77 & 6.54 & 4.94 & 3 \\
\textbf{LSR-LGBMW1} & 4.96 & 15.34 & 27.00 & 16.69 & 14.76 & 0 \\
grb-f1 & 5.21 & 28.24 & 50.47 & 28.36 & 22.72 & 0 \\
\bottomrule
\end{tabular}
\end{table}

\begin{table}[h]
\centering
\caption{\mkblue{Primal integral (\%) with ML models for 2-hour runs on dataset GT}}
\label{tab:primal_integral_ml_GT_2-hour}
\begin{tabular}{lccccccc}
\toprule
{Scenario} & \multicolumn{3}{c}{Quantiles} & {Mean} & {Geomean} & {Wins} \\
{} & {0.1} & {0.5} & {0.9} & {} & {} & {} \\
\midrule
LB-0.8 & 6.69 & 28.83 & 39.43 & 24.95 & 20.72 & 0 \\
LS* & 0.05 & 30.94 & 70.43 & 34.09 & 16.57 & 5 \\
LS*-3 & 0.26 & 5.13 & 16.38 & 6.93 & 5.50 & 12 \\
\textbf{LSR-LGBM} & 3.39 & 10.13 & 21.05 & 11.22 & 10.27 & 0 \\
\textbf{LSR-LGBMW-1} & 1.57 & 7.32 & 13.86 & 8.61 & 7.13 & 7 \\
\textbf{LSR-LGBMW1} & 5.75 & 16.80 & 25.88 & 16.73 & 14.99 & 0 \\
grb-f1 & 2.97 & 14.01 & 24.66 & 14.07 & 12.32 & 0 \\
\bottomrule
\end{tabular}
\end{table}

\begin{table}[h]
\centering
\caption{\mkblue{Shifted geomean of primal gap over time on dataset GT}}
\label{tab:fcn:gap_over_time_GT}
\begin{tabular}{lccccc}
\toprule
{Scenario} & \multicolumn{4}{c}{Time after probing} \\
{} & {1800s} & {2700s} & {3600s} & {4800s} \\
\midrule
LB-0.8 & 24.28 & 23.41 & 23.19 & 22.20 \\
LS* & - & 2.64 & 2.40 & 1.99 \\
LS*-3 & 8.45 & 8.30 & 8.30 & 8.30 \\
\textbf{LSR-LGBM} & 10.07 & 9.68 & 9.52 & 9.40 \\
\textbf{LSR-LGBMW-1} & 7.36 & 6.22 & 5.50 & 5.13 \\
\textbf{LSR-LGBMW1} & 15.66 & 14.93 & 14.81 & 14.77 \\
grb-f1 & 34.47 & 33.11 & 31.59 & 28.97 \\
\bottomrule
\end{tabular}
\end{table}

\begin{table}[h]
\centering
\caption{\mkblue{Primal gap (\%) with ML models for 2-hour runs on dataset B}}
\label{tab:primal_gap_ml_B_2-hour}
\begin{tabular}{lccccccc}
\toprule
{Scenario} & \multicolumn{3}{c}{Quantiles} & {Mean} & {Geomean} & {Wins} \\
{} & {0.1} & {0.5} & {0.9} & {} & {} & {} \\
\midrule
LB-0.8 & 0.07 & 0.33 & 0.56 & 0.32 & 1.31 & 1 \\
LS* & 0.00 & 0.04 & 0.28 & 0.10 & 1.10 & 11 \\
LS*-3 & 0.71 & 1.45 & 1.98 & 1.35 & 2.30 & 0 \\
\textbf{LSR-LGBM} & 0.04 & 0.14 & 0.41 & 0.20 & 1.19 & 2 \\
\textbf{LSR-LGBMW-1} & 0.00 & 0.01 & 0.19 & 0.06 & 1.06 & 14 \\
\textbf{LSR-LGBMW1} & 0.06 & 0.22 & 0.53 & 0.27 & 1.25 & 2 \\
grb-f1 & 0.08 & 0.40 & 0.67 & 0.40 & 1.38 & 1 \\
\bottomrule
\end{tabular}
\end{table}

\begin{table}[h]
\centering
\caption{\mkblue{Primal integral (\%) with ML models for 2-hour runs on dataset B}}
\label{tab:primal_integral_ml_B_2-hour}
\begin{tabular}{lccccccc}
\toprule
{Scenario} & \multicolumn{3}{c}{Quantiles} & {Mean} & {Geomean} & {Wins} \\
{} & {0.1} & {0.5} & {0.9} & {} & {} & {} \\
\midrule
LB-0.8 & 0.25 & 0.39 & 0.60 & 0.41 & 1.41 & 0 \\
LS* & 0.11 & 0.22 & 0.37 & 0.23 & 1.22 & 5 \\
LS*-3 & 0.75 & 1.22 & 1.45 & 1.19 & 2.17 & 0 \\
\textbf{LSR-LGBM} & 0.11 & 0.19 & 0.42 & 0.24 & 1.23 & 6 \\
\textbf{LSR-LGBMW-1} & 0.05 & 0.14 & 0.28 & 0.15 & 1.14 & 17 \\
\textbf{LSR-LGBMW1} & 0.12 & 0.23 & 0.53 & 0.29 & 1.28 & 4 \\
grb-f1 & 0.41 & 0.67 & 0.80 & 0.64 & 1.63 & 0 \\
\bottomrule
\end{tabular}
\end{table}

\begin{table}[h]
\centering
\caption{\mkblue{Shifted geomean of primal gap over time on dataset B}}
\label{tab:fcn:gap_over_time_B}
\begin{tabular}{lccccc}
\toprule
{Scenario} & \multicolumn{4}{c}{Time after probing} \\
{} & {1800s} & {2700s} & {3600s} & {4800s} \\
\midrule
LB-0.8 & 1.52 & 1.47 & 1.45 & 1.42 \\
LS* & 1.30 & 1.20 & 1.15 & 1.13 \\
LS*-3 & 2.41 & 2.37 & 2.33 & 2.32 \\
\textbf{LSR-LGBM} & 1.26 & 1.23 & 1.21 & 1.20 \\
\textbf{LSR-LGBMW-1} & 1.20 & 1.12 & 1.10 & 1.07 \\
\textbf{LSR-LGBMW1} & 1.31 & 1.28 & 1.27 & 1.26 \\
grb-f1 & 1.88 & 1.82 & 1.78 & 1.71 \\
\bottomrule
\end{tabular}
\end{table}

\begin{table}[h]
\centering
\caption{\mkblue{Primal gap (\%) with ML models for 2-hour runs on dataset D}}
\label{tab:primal_gap_ml_D_2-hour}
\begin{tabular}{lccccccc}
\toprule
{Scenario} & \multicolumn{3}{c}{Quantiles} & {Mean} & {Geomean} & {Wins} \\
{} & {0.1} & {0.5} & {0.9} & {} & {} & {} \\
\midrule
LB-0.8 & 0.28 & 0.46 & 0.69 & 0.49 & 1.48 & 0 \\
LS* & 0.00 & 0.12 & 0.34 & 0.14 & 1.13 & 6 \\
LS*-3 & 1.73 & 2.69 & 3.81 & 2.75 & 3.66 & 0 \\
\textbf{LSR-LGBM} & 0.00 & 0.16 & 0.36 & 0.18 & 1.17 & 4 \\
\textbf{LSR-LGBMW-1} & 0.00 & 0.00 & 0.41 & 0.11 & 1.10 & 16 \\
\textbf{LSR-LGBMW1} & 0.00 & 0.20 & 0.58 & 0.24 & 1.23 & 4 \\
grb-f1 & 0.48 & 0.99 & 1.26 & 0.90 & 1.88 & 0 \\
\bottomrule
\end{tabular}
\end{table}

\begin{table}[h]
\centering
\caption{\mkblue{Primal integral (\%) with ML models for 2-hour runs on dataset D}}
\label{tab:primal_integral_ml_D_2-hour}
\begin{tabular}{lccccccc}
\toprule
{Scenario} & \multicolumn{3}{c}{Quantiles} & {Mean} & {Geomean} & {Wins} \\
{} & {0.1} & {0.5} & {0.9} & {} & {} & {} \\
\midrule
LB-0.8 & 0.28 & 0.44 & 0.66 & 0.48 & 1.47 & 0 \\
LS* & 0.18 & 0.34 & 0.46 & 0.32 & 1.31 & 3 \\
LS*-3 & 1.63 & 2.35 & 2.96 & 2.37 & 3.32 & 0 \\
\textbf{LSR-LGBM} & 0.11 & 0.22 & 0.38 & 0.25 & 1.24 & 5 \\
\textbf{LSR-LGBMW-1} & 0.07 & 0.16 & 0.40 & 0.20 & 1.19 & 15 \\
\textbf{LSR-LGBMW1} & 0.14 & 0.24 & 0.57 & 0.30 & 1.29 & 7 \\
grb-f1 & 0.50 & 0.95 & 1.08 & 0.87 & 1.85 & 0 \\
\bottomrule
\end{tabular}
\end{table}

\begin{table}[h]
\centering
\caption{\mkblue{Shifted geomean of primal gap over time on dataset D}}
\label{tab:fcn:gap_over_time_D}
\begin{tabular}{lccccc}
\toprule
{Scenario} & \multicolumn{4}{c}{Time after probing} \\
{} & {1800s} & {2700s} & {3600s} & {4800s} \\
\midrule
LB-0.8 & 1.48 & 1.48 & 1.48 & 1.48 \\
LS* & - & 1.28 & 1.22 & 1.18 \\
LS*-3 & 3.74 & 3.71 & 3.69 & 3.68 \\
\textbf{LSR-LGBM} & 1.30 & 1.24 & 1.22 & 1.19 \\
\textbf{LSR-LGBMW-1} & 1.25 & 1.19 & 1.15 & 1.13 \\
\textbf{LSR-LGBMW1} & 1.34 & 1.29 & 1.27 & 1.25 \\
grb-f1 & 2.07 & 2.06 & 2.05 & 2.02 \\
\bottomrule
\end{tabular}
\end{table}

\section{Perspectives and Future Work}
\label{sec:conclusion}

The purpose of this study was to explore how to use an ML prediction to improve a SOA heuristic for the \fcn. To do so, we propose a methodology that takes advantage of supervised learning and off-the-shelf ML models. Given that the tools we use are highly accessible, we believe that our method is reasonably easy to reproduce and adapt to other problems similar to \fcn. Also, given that predictions are made on a local level, the method is not computationally expensive even for large instances. It constitutes a stepping stone for hybrid ML-OR approaches that still require further research to be fully operational. 

A noteworthy contribution of this research is the development of sampling routines that generate solutions for the \fcn\ in an efficient manner. Sampling routines are a key component of our method, as they provide the necessary information for the training process. We demonstrated that sampled solutions form highly informative features, which greatly increase the predictive power of the ML model. Furthermore, our study on feature sets revealed that graph features are not as useful as sampling features. ML models trained on graph features tend to overfit, and exhibit poor generalization. This is a key discovery that suggests that GNNs might not be necessary to solve CO problems even when the problem is defined on a graph.

Our experimental analysis contains important insights on the integration of ML predictions in a heuristic. We tested different ways to integrate the ML predictions in LS algorithms, and we showed that the best method is to use them to reduce the problem size by removing arcs from the graph. We found that it is most resilient to the noise in the ML predictions. This method is particularly effective on instances sampled from a uniform distribution, where it outperforms the original heuristic. We also showed that the ML model can be biased to improve the quality of the solutions. 

Two limitations of our ML-based approach are worth mentioning. First, the method does not outperform the SOA heuristic on GT instances. A human-designed heuristic remains the best option to find high-quality solutions to the \fcn\ when the computational time budget is large. A deterministic heuristic is also more reliable than our hybrid approach. ML predictions always come with a level of uncertainty, which translates into a risk of compromising the quality of the final solution. For this reason, we believe that the best use case for our ML approach is to help find high-quality solutions early in the solution process. This is useful within a hyperparameter tuning framework~\citep{hoos_automated_2012, su_automatic_2021} to find good configurations of a heuristic. Second, in the absence of direct indications from the sampled data, the ML model will fail to accurately predict the target.  Often, arcs appear to be closed given the sampled solutions, but they are actually open in the best solution. This is a consequence of the fact that sampling routines cannot definitively identify the optimal solution. This constitutes an opportunity for further research. One possible avenue is to replace the human-designed sampling routines by a ML model that has learned to explore the solution space more efficiently.

%

\newpage

\clearpage

\newpage

\bibliographystyle{apalike}

\bibliography{fcn_bib}

\end{document}